\documentclass[a4paper,12pt,legno]{article}

\usepackage{amsmath}
\usepackage{amsfonts}
\usepackage{amsthm}
\usepackage{mathrsfs}
\usepackage{amssymb}

\newtheorem{thm}{THEOREM}

\newtheorem{prop}{PROPOSITION}
\newtheorem{cor}{COROLLARY}

\theoremstyle{remark}

\usepackage{graphicx}
\usepackage[T1]{fontenc}
\usepackage[polish]{babel}
\usepackage[cp1250]{inputenc}

\newcommand{\Rset}{\mathbb{R}}

\newcommand{\F}{\mathcal{F}}

\begin{document}

\title{The strong mixing and the selfdecomposability properties}

\author{Richard C. Bradley\footnote{Department of Mathematics, Indiana University, Bloomington, Indiana 47405, USA.} \ and \
Zbigniew J. Jurek\footnote{Institute of Mathematics, University of
Wrocław, Pl. Grunwaldzki 2/4, 50-384 Wrocław, Poland.}\ \
\footnote{Research funded by Narodowe Centrum Nauki (NCN) grant no
\newline Dec2011/01/B/ST1/01257.}}

\date{July 21, 2013}

\maketitle

\maketitle

\begin{quote}\textbf{Abstract.}
It is proved that infinitesimal triangular arrays
obtained from normalized partial sums of strongly
mixing (but not necessarily stationary) random
sequences, can produce as limits only
selfdecomposable distributions.



\emph{Abbreviated title:} Strong mixing and selfdecomposability

\end{quote}

   Selfdecomposable probability measures
(in other words, the L\'evy class L distributions)
form (by definition) the class of possible limiting distributions of normalized partial sums from sequences
of independent (but not necessarily identically
distributed) random variables, under certain natural
technical assumptions on the normalizing constants.
The aim of this note is to show that selfdecomposable probability measures also form the class of possible
limiting distributions of normalized partial sums
from (not necessary stationary) strongly mixing
sequences, under the same technical assumptions on the normalizing constants.
The proof will utilize the standard Bernstein blocking
technique.

For normalized partial sums from \emph{strictly stationary},
strongly mixing random sequences, with a mild natural assumption on
the normalizing constants, two other classes of distributions ---
the stable and infinitely divisible laws --- have long been known to
play the same roles respectively as they do in the case of i.i.d.\
sequences: as possible limit laws (i) along the entire sequence of
normalized partial sums, and (ii) along subsequences of normalized
partial sums. For further information and references on those
classic results, see e.g.\ Volume 1, Chapter 12 of Bradley (2007).
We shall not treat the particular case of strict stationarity
further here.

\medskip
\medskip
\textbf{1. Notations and basic notions.}

Let $(\Omega,\F, P)$ be a probability space. Let $E$ be a real
separable Banach space, with norm $||\cdot||$ and Borel
sigma-algebra $\mathcal{E}$. By $\mathcal{P}\equiv\mathcal{P}(E)$ we
denote the set of all Borel probability measures on $E$, with
\emph{the convolution operation} denoted by $``\ast"$ and \emph{weak
convergence} denoted by $``\Rightarrow"$, which make $\mathcal{P}$ a
topological convolution semigroup.

Measurable functions $\xi: \Omega \to E$ are called \emph{Banach
space valued random variables}
(in short: $E$-valued rv's) and
$\mathcal{L}(\xi)(A):=P\{\omega\in\Omega: \xi(\omega)\in A\}$, for
$A\in\mathcal{E}$, is \emph{the probability distribution of $\xi$}.
Then for stochastically independent
$E$-valued random variables
$\xi_1$ and $\xi_2$ we have that $\mathcal{L}(\xi_1)\ast
\mathcal{L}(\xi_2)=\mathcal{L}(\xi_1+\xi_2)$. Also for $c\in
\Rset\setminus{\{0\}}$ and rv $\xi$ we define
$\mathcal{L}(c\,\xi)(A)=\mathcal{L}(\xi)(c^{-1}\,A)=:T_c(\mathcal{L}(\xi))(A)$,
for $A\in\mathcal{E}$. Similarly, for $T_c:E\to E$ given by
$T_cx:=cx$, we define $T_c\mu$, for $\mu\in\mathcal{P}$.
Hence $T_c(\mu\ast\nu)=T_c\mu \ast T_c\nu$.

For two sub-$\sigma$-fields
$\mathcal{A}$ and $\mathcal{B}$  of  $\F$ we define \emph{the measure
of dependence} $\alpha$ between them as follows:
\[
\alpha(\mathcal{A},\,\mathcal{B}):=\sup_{A\in\mathcal{A}, B\in
\mathcal{B}}\,|P(A\cap B)-P(A)P(B)|\,.
\]
For a given sequence $\textbf{X}:=(X_1,X_2,\dots)$ of
$E$-valued random variables, we define for each positive
integer $n$ the dependence coefficient
\begin{equation}
\alpha(n)\equiv \alpha(\textbf{X};n):=\sup_{j\in {\bf N}}\,\alpha\big(\,\sigma(X_k, 1\le k \le j), \, \,
\sigma(X_k, k \ge j+n)\big),
\end{equation}
where $\sigma(\dots)$ denotes the $\sigma$-field generated
by (\dots).
We will say that a sequence \textbf{X} is \emph{strongly mixing} (Rosenblatt (1956)) if
\begin{equation}
\alpha(n)\to 0 \ \ \mbox{as} \ \ n\to \infty.
\end{equation}
Of course, if
the elements of
\textbf{X} are stochastically independent then
$\alpha(\textbf{X};n)\equiv 0$.
Many known stochastic processes (including many
Markov chains, many Gaussian sequences, and many
models from time series analysis) have long been known
to be strongly mixing; see e.g.\ Bradley (2007).

\medskip
Suppose that for stochastically independent $E$-valued rv's $\xi_j,
j \in {\bf N}$, there exist sequences of real numbers $a_n$ and
vectors $b_n\in E$ and a probability measure $\nu$ such that
\begin{equation}
(i) \ \ a_n>0 \ \   \mbox{and} \  \ \forall(\epsilon>0) \
\lim_{n\to\infty}\max_{1\le k\le
n}P(\{\omega\in\Omega:\,a_n||\xi_k(\omega)||>\epsilon\})=0
\end{equation}
(the so called \emph{infinitesimality condition}) and
\begin{equation}
(ii)\ \ \lim_{n \to
\infty}\,P(\{\omega\in\Omega:\,
a_n(\xi_1+\xi_2+\dots+\xi_n)(\omega)+b_n
\in B\})=\nu(B)
\end{equation}
for every Borel set $B \subset E$ whose boundary
$\partial B$ satisfies $\nu(\partial B)=0$;
then the measure
$\nu$ is called \emph{selfdecomposable} or a
\emph{L\'evy class L
distribution}.

There are two basic characterizations of the class $L$:
\emph{the convolution decomposition} and
\emph{the random integral representation}.
The first one says that
\begin{equation}
[\ \nu \in L \ ]\  \ \ \mbox{iff} \ \ [\ \forall(0<c<1)\,\
\exists\,(\nu_c \,\in \mathcal{P}(E)) \ \ \nu=T_c\nu\ast\nu_c \ ],
\end{equation}
and hence the term  \emph{selfdecomposability}; cf.\ Jurek
and Mason (1993), Theorem 3.9.2.

The second one says that
\begin{equation}
[\ \nu \in L \ ]  \ \mbox{iff} \  \
\nu=\mathcal{L}\big(\int_0^{\infty}e^{-t}dY_{\rho}(t)\big),
\end{equation}
for some L\'evy process $(Y_{\rho}(t), t\ge 0)$  such that
$\mathcal{L}(Y_{\rho}(1))=\rho$ and  the log-moment $\mathbb{E}[\log
(1+||Y_{\rho}(1)||)]<\infty$; cf.\ Jurek and Vervaat (1983) or Jurek
and Mason (1993), Theorem 3.9.3. The L\'evy process $Y_{\rho}$ in
(6) is referred to as \emph{the background driving L\'evy process}
(in short: BDLP) of the selfdecomposable probability measures $\nu$.

Finally let us note that in terms of probability measures, (4) means
that
\[
 T_{a_n}(\rho_1\ast\rho_2\ast...\ast
\rho_n)\ast\delta_{b_n}\Rightarrow \nu \ \ \mbox{as} \ n\to\infty,
\]
where $\rho_i=\mathcal{L}(\xi_i)$ for $i=1,2,\dots$. For probability
theory on Banach spaces we refer to Araujo and Gin\'e (1980).

\medskip
\medskip
\textbf{2. Strong mixing and selfdecomposability.}

Here is the main result of this note:
\begin{thm}
Let $\textbf{X}:=(X_1,X_2,\dots)$ be a sequence of Banach space E
valued random variables with partial sums $S_n:=X_1+X_2+\dots+X_n$,
and let $(a_n)$ and $(b_n)$ be sequences of real numbers and
elements in $E$ respectively, and suppose the following conditions
are satisfied:

(i) $\alpha(n)\to 0$ as $n\to \infty$,  i.e.\ the sequence
$\textbf{X}$ is strongly
mixing;

(ii) $a_n>0$ and  the triangular array $(a_n\,X_j, \ 1 \le j \le n,
n\ge 1)$ is

infinitesimal;

(iii) $a_nS_n+b_n\Rightarrow \mu$ as $n \to \infty$
for some non-degenerate probability

measure $\mu \in \mathcal{P}$.

\noindent
Then the limit distribution $\mu$ is selfdecomposable.
\end{thm}
\emph{Proof of theorem.}  Our aim is to show that $\mu$ satisfies
the convolution decomposition (5). The argument below is divided
into a few steps/observations, some of which are quite elementary.

\medskip
\emph{Step 1.}   $a_n\to 0$ and ${a_{n+1}}/{a_n} \to 1$ as $n\to
\infty$.

\emph{Proof}. a) Let us first show that $a_n \to 0$.
Suppose instead that
$Q$ is an infinite subset of ${\bf N}$, $d > 0$, and
$a_n > d$ for all $n \in Q$.
Then for all $\epsilon >0$
\[
\max_{1\le k \le n} P(d|X_k|>\epsilon) \le \max_{1 \le k \le n} \ P(a_n|X_k|>\epsilon)\to 0 \ \mbox{as} \ n \to \infty,
n \in Q.
\]
Thus for each fixed $k$ and each $\epsilon>0$ we have
$P(|X_k|>\epsilon/d)=0$; i.e.\
for each $k \in {\bf N}$, $X_k=0$ with probability one.
So, $\mu$ is degenerate which contradicts the assumption (iii).
Thus $a_n \to 0$ after all.

b)  Now let us show that $a_{n+1}/a_n \to 1$.
Let $T_n:= a_nS_n+b_n\Rightarrow \mu$. Since, by the
infinitesimality assumption (ii), $a_{n+1}X_{n+1}\to 0$ in
probability, therefore (letting ``$\lim$'' denote
limit in distribution)
\[
\mu=\lim_{n\to\infty} T_{n+1}= \lim_{n\to \infty}
(a_{n+1}S_n+b_{n+1})=\lim_{n \to \infty}
\big[\frac{a_{n+1}}{a_n}\,\, T_n+
(b_{n+1}-b_n\frac{a_{n+1}}{a_n})\big]\,.
\]
Hence the Convergence of Types Theorem (see e.g.\
Proposition 2.7.1 in Jurek and Mason (1993))
gives that ${a_{n+1}}/{a_n} \to 1$ and
$b_{n+1}-b_n{a_{n+1}}/{a_n}\to 0$ as $n\to \infty$, which
completes the argument for Step 1.

\medskip
\emph{Step 2.} Now for what follows, let $c \in (0,1)$
be arbitrary but fixed.

For the sequence ($a_n$) let us
define
\[
m_n:=\max_{1\le k \le n-1}\{k: \frac{a_n}{a_k}\le c\}, \ \ \mbox{if
such $k$ exists, or}\ \ m_n:=1 \ \ \mbox{otherwise}.
\]
Then we have
\[
m_n\to \infty, \ \  n-m_n\to \infty \ \ \ \mbox{and} \ \ \
\frac{a_n}{a_{m_n}}\to c   \ \ \ \mbox{as} \  \ n\to \infty.
\]
\emph{Proof.} This is so, because from the definition
of $m_n$ we get
\[
\frac{a_n}{a_{m_n}}\le c< \frac{a_{n}}{a_{m_n+1}} \ \ \ \mbox{or} \
\ \ 1 \le c\, \frac{a_{m_n}}{a_n} < \frac{a_{m_n}}{a_{m_n+1}}  \ \ \
\mbox{for all sufficiently large}\ n.
\]

\medskip
\emph{Step 3.} For an infinitesimal triangular array $\{a_nX_k:
1\le k \le n, n=1,2,\dots\}$  (i.e.\ satisfying  the
condition (3))
there exists a non-increasing sequence of positive numbers
$\delta_1,\delta_2,\dots$, each $\leq 1$, such that
\[
(i)\ \ \delta_n \to 0 \  \ \mbox{as} \ \  n \to \infty \ \ \mbox{and}
\]
\[
(ii)\ \ \forall(n\ge 1)\ \ \forall (1\le k \le n) \ \
P(a_n\|X_k\| \ge \delta_n)\le \delta_n.
\]
\emph{Proof.} Simply note that (3) is equivalent to the following:
$\|X_{n,j_n}\| \to 0$ in probability, for all
$j_n\in {\bf N}$
such that $1\le j_n \le n$.

\medskip
\emph{Step 4.} For the sequences $(m_n)$ from Step 2, and
$(\delta_n)$ from Step 3, let us choose positive integers $q_n$ such
that
\[
 q_n \le \delta^{-1/2}_n \ \ \mbox{for all}\ n,\ \ q_n\to \infty\ \mbox{as}\ n \to \infty,\ \ \mbox{and}
\]
\[
 \  \ m_n+q_n<n  \ \ \mbox{for all sufficiently large} \ n\ge 1.
\]
Consequently, the partial sums $S_n$ can be written in three blocks
\begin{equation}
S_{m_n}+ (S_{m_n+q_n} - S_{m_n})+ (S_n -S_{m_n+q_n})=S_n \ \
\mbox{for sufficiently large $n\ge 1$}.
\end{equation}

\medskip
\emph{Step 5.} The
rv's $V_n:=a_n(S_{m_n+q_n} - S_{m_n})$ satisfy
$\|V_n\| \to 0$
in probability as $n \to \infty.$

\medskip
\emph{Proof.} This is so, because for any $\epsilon>0$, applying Steps 3 and 4,
for all $n$ such that $\delta_n^{1/2}< \epsilon$
and $m_n + q_n < n$,  we
have that
$\delta_n\,q_n \le \delta_n^{1/2}<\epsilon$
and thus
\begin{multline*}
P(||V_n||\ge \epsilon)=P(a_n||\sum_{k=m_n+1}^{m_n+q_n}X_k||\ge
\epsilon)
\le \sum_{k=m_n+1}^{m_n+q_n}P(a_n||X_k||\ge \epsilon/q_n)\\
\le \sum_{k=m_n+1}^{m_n+q_n}P(a_n||X_k||\ge
\epsilon\,\delta_n^{1/2})\le \sum_{k=m_n+1}^{m_n+q_n}P(a_n||X_k||\ge
\delta_n) \le \delta_n\, q_n \le \epsilon.
\end{multline*}

\medskip
\emph{Step 6.} For the sequence
$U_n:=({a_n}/{a_{m_n}})(a_{m_n}S_{m_n}+b_{m_n})$ we have that
\[
U_n\Rightarrow T_c\mu \ \ \ \  \mbox{and}\ \  \ U_n+V_n\Rightarrow
T_c\mu \,, \ \mbox{as} \ n \to \infty.
\]

\medskip
\emph{Proof.} First of all, by the assumption (iii) in the Theorem
1, and by Step 2 we have that $U_n\Rightarrow T_c\mu$. This with Step 5
(and e.g.\ Theorem 3.1 of Billingsley (1999))
proves the above claim.

\medskip
\emph{Step 7.} The family of probability distributions of
the random variables
\begin{equation*}
W_n:= a_n(S_n -S_{m_n+q_n})+b_n-\frac{a_n}{a_{m_n}}\,b_{m_n} , \  \
n\ge 1,
\end{equation*}
its tight. Equivalently, it is conditionally compact by the
Prokhorov Theorem.

\medskip
\emph{Proof.} Note that using (7) we have that
\begin{equation}
U_n+V_n +W_n=aS_n+b_n\Rightarrow \mu \ \mbox{as}  \ n\to \infty.
\end{equation}
By Prokhorov Theorem (see e.g.\ Theorem 1.6.6 in Jurek and Mason
(1993) or Theorem 2.10 in Araujo and Gin\'e (1980)),  using Step 6
and the formula (8), for a given $\epsilon>0$ there exist compact
sets $K_1$ and $K_2$ such that
\[
P(-(U_n+V_n)\in K_1)>1-\epsilon/2 \  \mbox{and} \ P(a_nS_n+b_n\in
K_2)>1-\epsilon/2
\]
for all $n\ge 1$. Since $W_n=(a_nS_n+b_n) +(-(U_n+V_n))$ and
\[
\{a_nS_n+b_n\in K_2 \ \mbox{and} \ \  -(U_n+V_n)\in K_1\}\subset
\{W_n\in K_1+K_2\}
\]
we conclude that
\[
P(W_n\notin K_1+K_2)\le P(-(U_n+V_n)\notin K_1)
+ P(a_nS_n+b_n \notin K_2) \le \epsilon
\]
for all $n\ge 1$.
Of course the set $K_1 + K_2$ is compact by a
standard argument.
Thus the family of probability distributions of
$(W_n)$ is tight (equivalently, conditionally compact in
the weak convergence topology).

\medskip
\emph{Step 8.} The probability measure $\mu$ (in Theorem 1) satisfies the convolution equation
$\mu=T_c\mu\ast\nu$ for some probability measure
$\nu\in\mathcal{P}$.

\medskip
\emph{Proof.} In view of Step 7, there exists
$\nu\in\mathcal{P}$
and a subsequence $Q\subset {\bf N}$ such that
\begin{equation}
W_n\Rightarrow \nu,  \ \ \ \mbox{as}  \ \ \  n\to \infty \ \ \mbox{and} \ \ n\in Q .
\end{equation}
From (8) and Steps 5 and 6 we conclude that
\begin{equation}
U_n+W_n\Rightarrow \mu,  \ \ \mbox{as}  \  n\to \infty.
\end{equation}
Since $q_n\to \infty$ therefore by the assumption (i) in Theorem 1,
we conclude
\begin{multline}
\alpha(\sigma(U_n), \sigma(W_n))\le \alpha(\sigma(X_k, 1\le k \le
m_n), \sigma(X_k, m_n+q_n+1\le k\le n)\\ \le \alpha(\textbf{X},
q_n+1)\to 0. \ \ \  \
\end{multline}
By (11), (9), Step 6, and Corollary 1 in the Appendix below
\begin{equation}
U_n+W_n\Rightarrow T_c\mu\ast\nu \ \ \mbox{as} \ \ n \to\infty,
\end{equation}
which with (10) completes the proof of the Step 8.

\medskip
Finally, since our argument can be repeated for each $0<c<1$
(note that the distribution $\nu$ in Step 8 depends on $c$),
we conclude by the convolution decomposition (5) that the measure
$\mu$ is indeed selfdecomposable, which completes the
proof of Theorem 1.

\vspace{0.25in} \textbf{3. Appendix.}

For ease of reference let us quote here
\begin{prop}
Let $X$ and $Z$ be two random elements on $(\Omega,\F,P)$ with
values in a separable Banach space $E$ and let $U$ be a uniformly
distributed on $(0,1)$ real valued random variable stochastically
independent of $X$ and $Z$.

Suppose further that there exist $\epsilon>0$, $\delta>0$,
a positive integer $N$, a Borel set $D \subset E$, and
points $a_1, a_2,..., a_N$ in $D$ such that
\[
P(X\in D)>1-\delta  \ \ \ \mbox{and} \ \ \forall ( a\in D) \,\exists
(1\le k \le N)\, ||a-a_k|| \le \epsilon.
\]

Then there exists an E-valued random variable $Y$ such that

(a) $Y$ is measurable with respect to $\sigma(X,Z,U)$,

(b) $Y$ is independent of $Z$,

(c) $Y\stackrel{d}{=}X$,

(d) $P(||X-Y|| > 2\epsilon)\le \delta + 4\,N^{1/2}\,\alpha
(\sigma(X),\,\sigma(Z))$.
\end{prop}
This is a corollary from Theorem 16.17 in Bradley (2007), volume 2,
p. 139, for  $S=E$ and $\mathcal{A}=\sigma(Z)$.

\begin{cor}
Suppose $(X_n)_n$ and $(Z_n)_n$ be two sequences of E-valued random
elements and $\mu$,$\nu$ be two probability Borel measures on $E$.

If $X_n\Rightarrow \mu$ and $Z_n\Rightarrow \nu$  as $n\to \infty$
and $ \lim_{n\to \infty}\alpha(\sigma(X_n),\sigma(Z_n))=0$ then

\medskip
(a) for every $\epsilon>0$ there exists positive integer
$n_{\epsilon}$
such that for every $n\ge n_{\epsilon}$ there exists
an $E$-valued random variable $Y_n$
that is independent of $Z_n$ and satisfies both
$Y_n\stackrel{d}{=}X_n$
and $P(||X_n-Y_n||>\epsilon)<\epsilon$.

\medskip
(b) random elements $(X_n, Z_n)\Rightarrow \mu\times \nu$ and
$X_n+Z_n\Rightarrow \mu\ast\nu$, as $n\to \infty$.
\end{cor}
\emph{Proof.} Applying the Prokhorov Theorem for the
sequence $(X_n)_n$
we have that for a given $\epsilon>0$ there exists a compact set $K\subset E$ such that
\[
P(X_n\in K)\, > 1-\epsilon/2 \ \ \mbox{for all} \ \ n\in
{\bf N}.
\]
From the compactness of K, there are finitely points $a_1,a_2,\dots, a_M \in K$ such that the
open balls $B(a_i,\epsilon/2)$, $i=1,2,\dots,M$
cover K and thus
for any $x\in K$ there is an $a_i$ such that
$||x-a_i||< \epsilon/2$.

Applying Proposition 1 for each pair $X_n, Z_n$ we get
a random variable $Y_n$ such
that $Y_n$ is independent of $Z_n$, $Y_n\stackrel{d}{=}X_n$ and
\[
P(||X_n-Y_n|| > \epsilon)\le \epsilon/2+ 4\,M^{1/2}\,\alpha
(\sigma(X_n),\,\sigma(Z_n)).
\]
Since $\lim_{n\to \infty}\alpha(\sigma(X_n),\sigma(Z_n))=0$, we conclude  that there exists a positive integer
$n_{\epsilon}$ such that $P(||X_n-Y_n||>\epsilon)<\epsilon$ for all
$n\ge n_{\epsilon}$. This completes the proof of part (a).

Part (b) follows from part (a) and Lemma 1.1 and Theorem 1.1 in
Chapter 3 in Parthasarathy (1967).

\end{document}